\documentclass[11pt]{article}
\usepackage{latexsym, amsfonts}
\newcommand{\be}{\begin{equation}}
\newcommand{\ee}{\end{equation}}
\newcommand{\bea}{\begin{eqnarray}}
\newcommand{\eea}{\end{eqnarray}}
\newcommand{\bean}{\begin{eqnarray*}}
\newcommand{\eean}{\end{eqnarray*}}
\newcommand{\brray}{\begin{array}}
\newcommand{\erray}{\end{array}}
\newcommand{\ben}{\begin{equation}{nonumber}}
\newcommand{\een}{\end{equation}{nonumber}}

\newtheorem{dfn}{Definition}[section]
\newtheorem{thm}[dfn]{Theorem}
\newtheorem{lmma}[dfn]{Lemma}
\newtheorem{ppsn}[dfn]{Proposition}
\newtheorem{crlre}[dfn]{Corollary}
\newtheorem{xmpl}[dfn]{Example}
\newtheorem{rmrk}[dfn]{Remark}
\newcommand{\bdfn}{\begin{dfn}}
\newcommand{\bthm}{\begin{thm}}
\newcommand{\blmma}{\begin{lmma}}
\newcommand{\bppsn}{\begin{ppsn}}
\newcommand{\bcrlre}{\begin{crlre}}
\newcommand{\bxmpl}{\begin{xmpl}}
\newcommand{\brmrk}{\begin{rmrk}}
\newcommand{\edfn}{\end{dfn}}
\newcommand{\ethm}{\end{thm}}
\newcommand{\elmma}{\end{lmma}}
\newcommand{\eppsn}{\end{ppsn}}
\newcommand{\ecrlre}{\end{crlre}}
\newcommand{\exmpl}{\end{xmpl}}
\newcommand{\ermrk}{\end{rmrk}}

\newcommand{\IC}{\mathbb{C}}

\newcommand{\IR}{\mathbb{R}}

\newcommand{\IT}{\mathbb{T}}

\newcommand{\cla}{{\cal A}}
\newcommand{\clb}{{\cal B}}
\newcommand{\clc}{{\cal C}}
\newcommand{\cld}{{\cal D}}
\newcommand{\cle}{{\cal E}}
\newcommand{\clf}{{\cal F}}

\newcommand{\clh}{{\cal H}}
\newcommand{\cli}{{\cal I}}

\newcommand{\cll}{{\cal L}}
\newcommand{\clm}{{\cal M}}

\newcommand{\clq}{{\cal Q}}

\newcommand{\cls}{{\cal S}}

\newcommand{\clw}{{\cal W}}

\def\a*{{\cal A}_{h,*}}
\def\B{{\cal B}(h)}
\def\B1{{\cal B}_1(h)}
\def\b{{\cal B}^{\rm s.a.}(h)}
\def\b1{{\cal B}^{\rm s.a.}_1(h)}

\newcommand{\ot}{\otimes}

\newcommand{\raro}{\rightarrow}

\newcommand{\lgl}{\langle}
\newcommand{\rgl}{\rangle}

\def \qed {$\Box$}

\begin{document}
\begin{center}
{\large {\bf Quantum symmetries and quantum isometries of compact metric spaces}}\\
by\\
{\large Debashish Goswami {\footnote{ Partially supported by a project on `Noncommutative Geometry and Quantum Groups' funded by Indian National Science Academy.}}}\\ 
{\large Stat-Math Unit, Kolkata Centre,}\\
{\large Indian Statistical Institute}\\
{\large 203, B. T. Road, Kolkata 700 108, India}\\
{e mail: goswamid@isical.ac.in}\\
\end{center}

\begin{abstract}
We prove that a compact quantum group with faithful Haar state which has a  faithful action on a compact  space must be a Kac algebra, with  bounded antipode and  the square of the antipode being identity. The main tool in proving this is the theory of ergodic quantum group action on $C^*$ algebras. 

Using the above fact, we also formulate a definition of isometric action of a compact quantum group on a compact metric space, generalizing the definition given by  Banica for finite metric spaces, and prove for certain special class of metric spaces the existence of the universal object in the category of those compact quantum groups which act isometrically and are `bigger' than the classical isometry group.
\end{abstract}

\section{Introduction}
 It is a very natural and interesting question to study quantum symmetries of classical spaces, particularly metric spaces. 
In fact, motivated by some suggestions of Alain Connes, S. Wang defined (and proved existence) of quantum group analogues of the classical 
symmetry or automorphism groups of various types of finite structures such as finite sets  
 and finite dimensional matrix algebras (see \cite{wang1}, \cite{wang2}), 
and then these quantum groups were investigated in depth by a number of mathematicians including Wang, Banica, Bichon and others 
(see, for example, \cite{ban1}, \cite{ban2}, \cite{bichon} and the references therein). However, it is important to extend these 
ideas and construction to the `continuous' or `geometric' set-up. In  a series of articles initiated by us in \cite{goswami} 
and then followed up in \cite{jyotish}, \cite{qorient}, \cite{qdisc} and other articles, we have formulated and studied quantum group analogues 
of the group of isometries (or orientation preserving isometries) of Riemannian manifolds, including in fact noncommutative geometric set-up in the sense 
of \cite{connes} as well. It remains to see whether such construction can be done in a metric space set-up, without assuming any finer geometric 
(e.g. Riemannian or spin) structures. This aim  is partially achieved in the present article, generalizing  Banica's formulation of quantum isometry groups of 
finite metric spaces. Indeed,  in \cite{qdisc}, we have proposed a natural definition of `isometric' action of a (compact) quantum group on an arbitrary 
compact metric space (extending Banica's definition which was given only for finite metric spaces), and showed in some explicit examples the existence 
of a universal object in the category of all such compact quantum groups acting isometrically on the given metric space. 
In the present paper we slightly modify this definition (for finite spaces it is still the same) and have been able to prove the existence of such a universal object for some special class of metric spaces. 

In fact, the first part of the paper does not go into the isometry condition and concentrates on the general aspects of quantum group actions on compact metric spaces. In this context, we have been able to prove some interesting results, which in particular imply that only Kac algebras can act faithfully on compact  spaces. 
After this, we assume that the compact space is metrizable and taking into account the metric $d$, give several equivalent formulations of the `quantum isomeric action'. However, we could not prove existence of a `quantum isometry group' in general, although in some special situations (which include typical Lie groups and homogeneous spaces, and also all finite groups endowed with invariant metric) we have been able to establish existence of such a quantum isometry group.

\section{Quantum groups and their actions}

A compact quantum group (CQG for short) is a  unital $C^*$ algebra $\cls$ with a coassociative coproduct 
(see \cite{woro1}, \cite{woro2}) $\Delta$ from $\cls$ to $\cls \ot \cls$ (injective tensor product) 
such that each of the linear spans of $\Delta(\cls)(\cls \ot 1)$ and that of $\Delta(\cls) (1 \ot \cls)$ are norm-dense in $\cls \ot \cls$. 
From this condition, one can obtain a canonical dense unital $\ast$-subalgebra $\cls_0$ of $\cls$ on which linear maps $\kappa$ and 
$\epsilon$ (called the antipode and the counit respectively) making the above subalgebra a Hopf $\ast$ algebra. In fact, we shall always choose this dense 
 Hopf $\ast$-algebra to be the algebra generated by the `matrix coefficients' of the (finite dimensional) irreducible unitary representations (to be defined 
 shortly) of the CQG. The antipode is an anti-homomorphism and also satisfies $\kappa(a^*)=(\kappa^{-1}(a))^*$ for $a \in \cls_0$.
 
 It is known that there is a unique state $h$ on a CQG $\cls$ which is bi invariant in the sense that $({\rm id} \ot h)\circ \Delta(a)=(h \ot {\rm id}) \circ \Delta(a)=h(a)1$ for all $a$. The Haar state need not be faithful in general, though it is always faithful on $\cls_0$ at least. One also has $h(\kappa(a))=h(a)$ for $a \in \cls_0$. We also recall from \cite{woro2} that there exists a canonical one-parameter family $f_z$ indexed by $z \in \IC$ of 
  multiplicative linear functionals on $\cls_0$, with interesting properties listed in Theorem 1.4 of \cite{woro2}. In particular,  $f_z$ forms a one parameter group with respect to the `convolution' $\ast$ (i.e. $f_{z+z^\prime}=f_z \ast f_{z^\prime}:=(f_z \ot f_{z^\prime})\circ \Delta$), and $\kappa^2(a)=f_{-1} \ast a \ast f_1$. The following fact, which is  contained in Theorem 1.5 of \cite{woro2}, will be quite useful for us:
 \bppsn
 \label{unimod}
 The Haar state $h$ is tracial, i.e. $h(ab)=h(ba)$ for all $a,b$, if and only if $k^2={\rm id}$, which is also equivalent to $f_z=\epsilon~\forall z$.  In such a case, $\cls$ is called  a Kac algebra.
 \eppsn

We also need the following:
\blmma
\label{1234}
If the tracial Haar state of a Kac algebra  $\cls$ is faithful,  then   its antipode $\kappa$ admits a norm-bounded extension on $\cls$ satisfying $\kappa^2={\rm id}$. 
\elmma
{\it Proof:}\\ It follows from Theorem 1.6, (4) of \cite{woro2}, by noting that $\tau_{i/2}={\rm id}$ for kac algebras, since $f_z=\epsilon$ for all $z$ in this case. \qed

\bcrlre
\label{kac}
Let $\cls$ be a CQG with faithful Haar state $h$ and assume that the mutiplicative functionals $f_z$ are identically equal to the counit $\epsilon$ on a norm-dense $\ast$-subalgebra $\cls_1$ of $\cls_0$ ($\cls_1$ may be strictly smaller than $\cls_0$). Then $\cls$ must be a Kac algebra.
\ecrlre
{\it Proof:}\\ 
In the notation of \cite{woro2}, we observe from the proof of Theorem 1.6 of \cite{woro2} that $\sigma_t(a)=a$ for all $a \in \cls_1$, where $\sigma_t$ denotes the $\ast$-automorphism defined on the whole of $\cls$ mentioned in (3) of Theorem 1.6 of \cite{woro2}. By the boundedness of $\sigma_t$ and the norm-density of $\cls_1$, we conclude that $\sigma_t={\rm id}$ on the whole of $\cls$, hence $h$ must be tracial. \qed

We say that a CQG $\cls$ (with a coproduct $\Delta$) (co)acts on a unital $C^*$ algebra $\clc$ if there is a unital $C^*$-homomorphism 
$\beta : \clc \raro \clc \ot \cls$ such that ${\rm Span}\{ \beta(\clc)(1 \ot \cls)\}$ is norm-dense in $\clc \ot \cls$, and it satisfies the coassociativity
 condition, i.e. $(\beta \ot {\rm id}) \circ \beta=({\rm id} \ot \Delta) \circ \beta$. It has been shown in \cite{podles} that 
 there is a unital dense $\ast$-subalgebra $\clc_0$ of $\clc$ such that $\beta$ maps $\clc_0$ into $\clc_0 \ot_{\rm alg} \cls_0$ (where $\cls_0$ is the dense 
 Hopf $\ast$-algebra mentioned before) and we also have 
   $({\rm id} \ot \epsilon)\circ
 \beta={\rm id}$ on $\clc_0$. In fact, this subalgebra $\clc_0$ comes from the canonical decomposition of $\clc$ into subspaces on each of 
which the action 
 $\beta$ is equivalent to an irreducible representation. More precisely, $\clc_0$ is the algebraic direct sum of finite dimensional vector spaces $\clc^\pi_i$, say, where $i$ runs over some index set $J_i$, and $\pi$ runs over some subset (say $T$) of the set of (inequivalent) irreducible unitary representations of $\cls$, and the restriction of $\beta$ to $\clc^\pi_i$ is equivalent to the representation $\pi$. Let $\{ a^{(\pi,i)}_j,j=1,...,d_\pi \}$ (where $d_\pi$ is the dimension of the representation $\pi$) be a basis of $\clc^\pi_i$ such that $\beta(a^{(\pi,i)}_j)=\sum_k a^{(\pi,i)}_k \ot t^\pi_{jk}$, 
  for elements $t^\pi_{jk}$ of $\cls_0$. The elements $\{ t^\pi_{jk}, \pi \in T;~j,k=1,...,d_\pi \}$ are called the `matrix coefficients' of the action $\beta$.
 
 We say that the action $\beta$ is faithful if the $\ast$-subalgebra of $\cls$ generated by elements of the form $(\omega \ot {\rm id})(\beta(a))$, where $a \in \clc$, $\omega$ being a bounded linear functional on $\clc$, is norm-dense in $\cls$. Since every bounded linear functional on a $C^*$ algebra is a linear combination of states, it is clear that in the definition of faithfulness, we can replace `bounded linear functional' by `state'. 
 \blmma
 Given an action $\beta$ of a CQG $\cls$ on $\clc$, with $\clc_0$, $\cls_0$ etc. as above, the following are equivalent:\\
 (i) The action $\beta$ is faithful.\\
 (ii) The $\ast$-algebra generated by the matrix coefficients is norm-dense in $\cls$.\\
 (iii) The $\ast$-algebra generated by elements of the form $(\omega \ot {\rm id})(\beta(a))$, where $a \in \clc_0$, and $\omega$ is a (not necessarily bounded) linear functional on $\clc_0$, is norm-dense in $\cls$.
 \elmma
 {\it Proof:}\\
 The equivalence of (ii) and (iii) is quite obvious, and so is the implication (i) $\Rightarrow$ (iii). To prove (ii) $\Rightarrow$ (i), fix any $\pi,i,j,k$ and consider (by Hahn-Banach theorem) a bounded linear functional $\omega$ on $\clc$ such that $\omega(a^{(\pi,i)}_k)=1$, $\omega(a^{(\pi,i)}_l)=0$ for $l \neq k$. Clearly, $(\omega \ot {\rm id})(\beta(a^{(\pi,i)}_j))=t^\pi_{jk}$. \qed\\


For a Hilbert $C^*$  module $\cle$ over a $C^*$ algebra $\clc$, we shall denote by $\cll(\cle)$ 
the $C^*$ algebra of adjointable $\clc$-linear maps from $\cle$ to $\cle$.    
We shall typically consider the Hilbert $C^*$  modules of the form $\clh \ot \cls$, where $\cls$ is a $C^*$ algebra
and the Hilbert module is the completion of $\clh \ot _{\rm alg} \cls$ w.r.t. the weakest topology which makes $\clh \ot_{\rm alg} \cls \ni X \mapsto
 \lgl X, X \rgl^{\frac{1}{2}} \in \cls$ continuous in norm. 
We shall use two kinds of  `leg-numbering' notation: for $T \in \cll(\clh \ot \cls)$, we denote by $T_{23}$ and $T_{13}$ the elements of 
 $\cll(\clh \ot \clh \ot \cls)$  given by  $T_{23}=I_\clh \ot T$, $T_{13}=\sigma_{12} \circ T_{23} \circ \sigma_{12}$, where $\sigma_{12}$ flips two copies of 
 $\clh$. On the other hand, we shall denote by $T^{12}$ and $T^{13}$ the elements $T \ot {\rm id}_\cls$ and $\sigma_{23} \circ T^{12} \circ \sigma_{23}$ respectively,
  of $\cll(\clh \ot \cls \ot \cls)$, where $\sigma_{23}$ flips two copies of $\cls$.

A unitary representation of  a CQG $(\cls, \Delta)$ in a Hilbert space $\clh$ is given by a unitary $U$ from $\clh$ to $\clh \ot \cls$, 
or equivalently, the unitary $\tilde{U} \in \cll(\clh \ot \cls)$ defined by $\tilde{U}(\xi \ot b)=U(\xi)(1 \ot b)$, 
for $\xi \in \clh, b \in \cls,$ satisfying $({\rm id} \ot \Delta) (\tilde{U})=\tilde{U}^{12}\tilde{U}^{13}$. 
 We denote by ${\rm ad}_U$ the map $\clb(\clh) \ni X \mapsto \tilde{U}(X \ot 1){\tilde{U}}^*$. 
We say that an action $\beta$ of  $\cls$ on a $C^*$ algebra $\clc$  is implemented by a unitary representation in a Hilbert space $\clh$ if there is a faithful  representation $\pi$ of $\clc$ in $\clb(\clh)$ such that $(\pi \ot {\rm id}) \circ \beta(a) ={\rm ad}_{U}(\pi(a))$ for all $a \in \clc$.  Given such a unitary representation $U$, 
we denote by $\tilde{U}^{(2)}$ the unitary $\tilde{U}_{13}\tilde{U}_{23} \in \cll(\clh \ot \clh \ot \cls)$, and consider ${{\rm ad}_U}^{(2)} \equiv {\rm ad}_{U^{(2)}}$, which is given by
 ${{\rm ad}_{U}}^{(2)}(x \ot y):={\tilde{U}}^{(2)}(x \ot y \ot 1){{\tilde{U}}^{(2)^*}}$, for $x, y \in \clb(\clh)$.  We use similar notation for 
an action $\beta$ of $\cls$ on some $C^*$  algebra $\clc$ implemented by $U$, i.e. take $\beta^{(2)} \equiv \beta^{(2)}_U$ (this may in general depend on $U$) to be the restriction of ${{\rm ad}_U}^{(2)}$ to $\clc \ot \clc$. This will be referred to as the `diagonal action', since in the commutative case, i.e.
 when $\clc=C(X)$, $\cls=C(G)$, with $G$ acting on $X$, the action $\beta^{(2)}$ does indeed correspond to the diagonal action of $G$ on $X \times X$. 
However, we warn the reader that when $\cls$ is no longer commutative as a $C^*$ algebra, i.e. not of the form $C(G)$ for some group $G$, $\beta^{(2)}$ may not leave $\clc \ot \clc$ (or even $\clc \ot_w \clc$) invariant, so may not be an action of $\cls$ on $\clc \ot \clc$. 

We shall denote by $\beta_{(2)} \equiv \beta_{(2),U}$ the $\ast$-homomorphism $ {\rm ad}_W$, where $W=\tilde{U}_{23}\tilde{U}_{13}$.
\brmrk
The `diagonal action' $\beta^{(2)}$ is not same as the one considered in \cite{qdisc}; in fact, the diagonal map of \cite{qdisc} is actually (at least for finite spaces) the unitary $U^{(2)}$ considered in the present paper, so is not an algebra homomorphism in general. 

\ermrk 

 s
\section{Necessity of the Kac algebra condition}
The aim of this section is to prove that if a CQG with faithful Haar state acts facithfully on $C(X)$, then the CQG must be Kac algebra. The main tool is the  results known about ergodic actions of CQG. Let is recall that an action $\beta$ of a CQG $\cls$ on a $C^*$ algebra $\clc$ is called ergodic if the fixed point subalgebra is trivial, i.e. $\beta(x)=x \ot 1$ for some $x\in \cla$ implies that $x$ is a scalar multiple of the identity element of $\cla$. We refer the reader to \cite{erg} for a detailed analysis of such ergodic actions. The following is an easy corollary of the results obtained in \cite{erg} (see also \cite{vaes}):
\bthm
\label{ergodic}
Let $\clc$ be a commutative unital $C^*$ algebra on which a CQG $\cls$ has an ergodic action $\beta$. Let $\clc_0$, $\cls_0$ be the dense $\ast$-subalgebras mentioned before, with the counit $\epsilon$ and antipode $\kappa$, and let $f_z$ be the one-parameter family of multiplicative maps discussed in Section 2. Then  $f_z(x)=\epsilon( x)$ for all $x$ belonging to the subalgebra of $\cls_0$  spanned by elements of the form $(\theta \ot {\rm id})(\beta(a))$, for $a \in \clc_0$ and $\theta$ being a linear functional on $\clc_0$. 
\ethm
{\it Proof:}\\
Let us recall the vector space decomposition of $\clc_0$ into subspaces $\clc^\pi_i$  discussed before,  and denote by $F_\pi$  (as in \cite{erg}, and also \cite{woro2}, where the symbol $F^\pi$ was used)   the operator 
 $({\rm id} \ot f_1) \circ \beta$ restricted to the span of $\clc^\pi_i$'s, for any fixed $\pi,i$. Note that the operator essentialy depends on $\pi$ only, as it is given on $a_j^{(\pi,i)}$ by $\sum_k a^{(\pi,i)}_k f_1(t^\pi_{j,k})$, so that we did not keep $i$ in the suffix. It is clear that  it suffices for us to prove that $F_\pi=I$ for each  $\pi$. It is in fact a positive invertible finite-dimensional operator satisfying ${\rm Tr}(F_\pi)={\rm Tr}(F_\pi^{_1})$ (see \cite{woro2}). Now by Proposition 18 of \cite{erg}, there is a unique faithful $\cls$-invariant state $\omega$ on $\clc$ and multiplicative linear map $\Theta: \clc_0 \raro \clc_0$ satisfying 
  $\omega(x \Theta(y))=\omega(yx)$ for all $x,y \in \clc_0$.  But $\clc$ being commutative, we have $\omega(x \Theta(y))=\omega(xy)$, and by the faithfulness of $\omega$, we must have $\Theta(y)=y$ for all $y \in \clc_0$. Moreover, by \cite{erg}, for each $\pi$ and $i$,  the restriction pf $\Theta$ to  $\clc^\pi_i$ coincides with a scalar multiple of $F_\pi$, and hence (noting also that ${\rm Tr}(F_\pi)={\rm Tr}(F_\pi^{-1})$) we get  $F_\pi=I$. \qed\\

We now prove the main result of this section.
\bthm
\label{unimod_main}
Let $(\cls, \Delta)$ be a CQG with faithful Haar state which acts faithfully on $C(X)$ for a compact space $X$. Then $\cls$ must be a Kac algebra.
\ethm
{\it Proof:}\\
Let $\kappa, \epsilon$ be the antipode and counit of $\cls$ (defined at least on $\cls_0$) respectively, and let $f_z$ be the one-parameter group of multiplicative functionals discussed before. 
For any $x \in X$, let $\beta_x$ denote the $\ast$-homomorphism $({\rm ev}_x \ot {\rm id}) \circ \beta$ from $C(X)$ to $\cls$, and let the range of $\beta_x$ be denoted by $\cls^x$.  Clearly, $\cls^x$ is a unital commutative $C^*$-algebra and by the coassociativity of $\beta$, it follows that $\Delta(\cls^x) \subseteq \cls^x \ot \cls$. Thus, $\Delta|_{\cls^x}$ gives an action of $\cls$ on $\cls^x$, and we claim that this is ergodic. Indeed, if $\Delta(a)=a \ot 1$ for $a \in \cls^x$, by appying the Haar state $h$ (say) of $\cls$ on the second copy of the tensor product we get $h(a)1=a$, i.e. $a$ is a scalar multiple of $1$. 

Now, let $\clc_0$ be the dense $\ast$-subalgebra of $C(X)$ on which $\beta$ is algebraic, and by the assumption of faithfulness, the $\ast$-algebra generated by $\cls^x_0:=\beta_x(\clc_0)$, with $x$ varying in $X$, is dense in $\cls$. Since the action $\Delta|_{\cls^x}$ is ergodic, we conclude  by Theorem \ref{ergodic} that $f_z(b)=\epsilon(b)$ $\forall z$, for all elements  $b$ of the form $b=(\omega \ot {\rm id})(\Delta(a))$, for any linear functional $\omega$ on $\cls^x_0$, and $a \in \cls^x_0$. Taking $\omega=\epsilon$,  (which is defined on $\cls^x_0 \subseteq \cls_0$), we see that $f_z=\epsilon$ on $\cls^x_0$ for every $x$, and using the facts that $f_z(ab)=f_z(a)f_z(b)$ and $f_z(a^*)=\overline{f_{-z}(a)}$ for all $a,b$, and $\epsilon$ is $\ast$-homomorphism, 
 we get that $f_z=\epsilon$ on  the $\ast$-subalgebra (say $\cls_1$) of $\cls_0$ generated by $\cls^x_0$'s, which is norm-dense in $\cls$ by faithfulenss. The theorem now follows from Corollary \ref{kac}.
\qed

\section{Quantum group of isometries of $(X,d)$}
\subsection{Definition of isometric action of compact quantum groups}
We have already noted that for any $C^*$ action $\beta$ of a CQG $\cls$  with faithful Haar state on $C(X)$, the antipode, say $\kappa$ of $\cls$ is defined and bounded on the $C^*$ subalgebra generated by $\beta(f)(x) \equiv ({\rm ev}_x \ot {\rm id})\circ \beta(f)$, $f \in C(X),$ $x \in X$. So $({\rm id} \ot \kappa)\circ \beta$ is a well-defined and norm-bounded map on $C(X)$. 

In view of this, it is natural to make the following definition:
\bdfn
\label{isodef}
Given an action $\beta$ of a CQG $\cls$ (with faithful Haar state) on  $\clc=C(X)$ (where $(X,d)$ is a compact metric space), we say that $\beta$ is  `isometric'  if the metric 
 $d \in C(X) \ot C(X)$ satisfies the following:
\be \label{isodef_formula} ({\rm id}_\clc \ot \beta)(d)=\sigma_{23} \circ (({\rm id}_\clc \ot \kappa) \circ \beta \ot {\rm id}_\clc)(d),\ee
 where $\sigma_{23}$ denotes the flip of the second and third tensor copies.
\edfn

\bthm
\label{isometry_dfn}
Given a $C^*$-action $\beta$ of a CQG $\cls$ (with faithful Haar state) on $C(X)$, the following are equivalent:\\
(i) The action is isometric.\\
(ii) $\forall x, y \in X$, one has $\beta(d_x)(y)=\kappa(\beta(d_y)(x))$, where $d_x(z):=d(x,z)$.\\
(iii) For some (hence all) unitary representation $U$ of $\cls$ on $\clh_U$ which implements $\beta$,  we have $\beta_U^{(2)}((\pi_U \ot \pi_U)(d))=(\pi_U \ot \pi_U)
(d) \ot 1$, where $\pi_U: \cls \raro \clb(\clh_U)$ denotes the imbedding of $\cls$ into $\clb(\clh_U)$.\\
(iv) For some (hence all)  unitary representation $U$ of $\cls$ on $\clh_U$ which implements $\beta$,  we have $\beta_{(2), U}((\pi_U \ot \pi_U)(d))=(\pi_U \ot \pi_U)
(d) \ot 1$.\\
\ethm
{\it Proof:}\\
The equivalence of (i) and (ii) is a consequence of  the continuity of the map $x \mapsto d_x \in C(X)$, and hence
 (by the norm-contractivity of $\beta$), the continuity of $x \mapsto \beta(d_x) \in C(X) \ot \cls$. Finally, to prove the equivalence of (i) and (iii), 
 we need to first note that in (i), i.e. the condition (\ref{isodef_formula}), $\beta$ can be replaced by 
 ${\rm ad}_U$ for any unitary representation $U$ of the CQG $\cls$ in some Hilbert space, which implements $\beta$. 
 Then, using the observation  that 
 ${(\tilde{U})}^{-1} (\pi_U(f) \ot 1)\tilde{U}=(\pi_U \ot \kappa)(\beta(f))$, we see that (i) is equivalent to 
 the following:
$$ \tilde{U}_{23} ((\pi_U \ot \pi_U)(d)\ot 1))\tilde{U}_{23}^{-1}=\tilde{U}_{13}^{-1}((\pi_U \ot \pi_U)(d)\ot 1)\tilde{U}_{13},$$ 
which is clearly 
 nothing but (iii).  
 
 Finally, the equivalence of (i) and (iv) follows from the symmetry of $d$, i.e. $d \sigma=d$, where $\sigma \in C(X \times X)$ is the map 
  $\sigma(x,y)=(y,x).$ To elaborate this calculation little more, let us denote by $\sigma_{12}$ the map $\sigma \ot {\rm id}_\cls$, and observe that 
    $\sigma_{12}\circ {\rm ad}_{U_{13}}\circ \sigma_{12}={\rm ad}_{U_{23}}$ and $\sigma_{12}\circ {\rm ad}_{U^{-1}_{23}}\circ \sigma_{12}=
    {\rm ad}_{U^{-1}_{13}}$. Now, it is clear that (iv) is obtained from (i) by replacing $d$ by $\sigma(d)$  in (\ref{isodef_formula}) and then also applying $\sigma_{12}$ on both sides of it. 
\qed\\

\brmrk
In case $\cls=C(G)$ and $\beta$ corresponds to a topological action of a compact group $G$ on $X$, it is clear that the condition (ii) above is nothing but the requirement $d(x, gy)=d(y,g^{-1}x) (=d(g^{-1}x,y)) $, which is obviously the usual definition of isometric group action.
\ermrk
\brmrk
For a finite metric space $(X,d)$, the present definition does coincide with Banica's definition in \cite{ban1} as well as the one proposed in \cite{qdisc}. Indeed, for such a space, $C(X)=l^2(X)$, and we can be choose $\clh_U=l^2(X)$ and $U$ to be the natural representation of $\cls$ coming from the action, with $d$ viewed both as an element of $C(X \times X)$ as well as of $l^2(X \times X)$. There is also the identically $1$ function, say ${\bf 1}$, in $l^2(X \times X)$, which is cyclic and separating  for $C(X \times X)$. Thus, the requirement $\beta^{(2)}_U(d)=d \ot 1$ is clearly equivalent to $U(d)=d \ot 1$, since ${U}^{-1} {\bf 1}={\bf 1} \ot 1$. 
This is precisely the proposed condition of \cite{qdisc} (and equivalent to the definition of Banica, as observed in \cite{qdisc}).
\ermrk
\brmrk
\label{imp_rem}
In the more general situation, assume that the CQG $\cls$ has faithful Haar state (otherwise one may replace the original CQG by a suitable quantum subgroup). Then consider any faithful state $\phi$ on $C(X)$ (given by integration w.r.t. some Borel probability measure $\mu$, say), and 
by averaging w.r.t. the Haar state, we get another faithful (and $\cls$-invariant) state, say $\overline{\phi}$, with the corresponding measure being $\overline{\mu}$. Clearly, the action $\beta$ extends to a unitary representation, say $U$, on $\clh:=L^2(X, \overline{\mu})$ which implements $\beta$. Moreover, since $\overline{\mu}$ is a Borel probability measure, we have $C(X) \subseteq L^2(X, \overline{\mu})$. and ${\bf 1}$ is a  cyclic separating vector for $C(X \times X)$ in  $L^2(X \times X)$ as before, such that ${U^{(2)}}^{-1}({\bf 1})={\bf 1} \ot 1$. Thus, $\beta$ is isometric in our sense if and only if $U^{(2)}(d)=d \ot 1$, $d$ being viewed as a vector in $L^2(X \times X)$. Similarly, using condition (v) of Theorem \ref{isometry_dfn}, we get $U_{23}U_{13}(d)=d \ot 1$. More generally, for any function $\phi(d)$, where $\phi$ is a continuous real-valued function on $\IR^+$, we have 
 $\beta^{(2)}(\phi(d))=\phi(d) \ot 1=\beta_{(2)}(\phi(d)),$ and hence also $U_{13}U_{23}(\phi(d))=U_{23}U_{13}(\phi(d))=\phi(d) \ot 1$.
\ermrk
\brmrk
In the recent article  \cite{sabbe}, the authors have generalized the notion of isometric action of CQG to the framework of Rieffel's compact quantum metric spaces. However, it is not yet clear whether their definition of isometric CQG-action is the same as the one given by us for a general compact metric space, although for finite spaces their equivalence has been proved by them. It will be also interesting to see whether analogues of the equivalent conditions (i)-(iv) proved above can be generalized too. 
\ermrk

\subsection{Existence of a universal isometric action}
It is a natural question to ask: does there exist a universal object in the category (say ${\bf Q}_{X,d}$) of all CQG acting isometrically (in our sense) on $(X,d)$? For finite metric spaces, the answer is clearly affirmative, and the universal object is the quantum isometry group defined by Banica. We are not yet able to settle this question in full generality. However, for a class of metric spaces which include typical Lie groups and their homogeneous spaces, we shall give an affirmative answer to a slightly modified question, namely, we shall prove the existence of a universal object in the subcategory of ${\bf Q}_{X,d}$ consisting of those for which  the classical isometry group $ISO(X,d)$ (viewed as a CQG) is a sub-object. This formulation of quantum isometry group tacitly assumes that such a quantum group, if exists,  should be `bigger' than the classical isometry group. 

\bthm
\label{exist}
Let us assume that there is a regular Borel probability measure $\mu$ on $X$ which is the unique $G$-invariant (where $G:=ISO(X)$ denotes the classical 
 isometry group) probability measure on $X$. Moreover, assume that there is a strictly positive, bijective function $\phi: \IR_+ \raro \IR_+$ such that  the integral operator $T$ given by the integral kernel $k(x,y)=\phi(d(x,y))$ in $L^2(X,\mu)$ is compact, with eigenvectors belonging to $C(X)$ and their linear span being norm-dense in $C(X)$. 
 
 Under the above assumptions,   there exists  a universal object, to be denoted by $QISO(X,d)$, in the subcategory of ${\bf Q}_{X,d}$ consisting of those for which  the classical isometry group $ISO(X,d)$ (viewed as a CQG) is a sub-object.
\ethm
{\it Proof:}\\
Let $\tau_\mu$ denote the state on $C(X)$ given by integration w.r.t. $\mu$.  Let $\clq$ be a CQG in the category mentioned in the statement of the theorem, with $\alpha$ being its action on $C(X)$ and $h$ being the Haar state on $\clq$. Since it contains $C(G)$ as a quantum subgroup, it is clear that the state $\tau_\mu \ast h:=(\tau_\mu \ot h) \circ \alpha$, which is $\clq$-invariant, will be $G$-invariant too. By the assumed uniqueness of $G$-invariant regular Borel probability measure, we conclude that $\tau_\mu=\tau_\mu \ast h$, 
 and now it is easy to argue, using the fact that $h$ is the Haar state of $\clq$, that $\tau_\mu$ is in fact $\clq$-invariant. Thus, the action $\alpha$ naturally 
  extends to a unitary representation on the Hilbert space  $\clh=L^2(\mu) $, which we denote by $U$, and $\tilde{U}$ denotes the corresponding unitary on the Hilbert module $\clh \ot \clq$.  
  
  We claim that \be \label{***}  U \circ T=(T \ot {\rm id}) \circ U. \ee  To this end, we first note that since $k=\phi(d)$, it follows from Remark \ref{imp_rem} that \be \label{1972} U_{23}U_{13}(k)=k \ot 1_\clq =U_{13}U_{23}(k).\ee  Now, note that $ Tf=({\rm id} \ot \tau_\mu)(k(1 \ot f))$ for $f \in C(X) \subseteq L^2(X, \mu)$. Thus, for $f,g \in C(X)$, we have the following (where $< \cdot, \cdot>_\clq$ will denote the $\clq$-valued inner product of the Hilbert modules $\clh \ot \clq$ or  $\clh \ot \clh \ot \clq$):
  \bean \lefteqn{<g \ot 1_\clq, U(Tf) >_\clq}\\
  &=& <g \ot 1 \ot 1_\clq, \widetilde{U_{13}}((1 \ot f)k \ot 1_\clq)>_\clq\\
  &=& <g \ot 1 \ot 1_\clq, (1 \ot f \ot 1_\clq) \widetilde{U_{13}}(k \ot 1_\clq)>_\clq\\
  &=& <\widetilde{U_{23}}(g \ot 1 \ot 1_\clq), \widetilde{U_{23}}(1 \ot f \ot 1_\clq){\widetilde{U_{23}}}^{-1}\widetilde{U_{23}}\widetilde{U_{13}}(k \ot 1_\clq) >_\clq\\
  &=& <g \ot 1 \ot 1_\clq, (1 \ot \alpha(f))\widetilde{U_{23}}\widetilde{U_{13}}(k \ot 1_\clq)>_\clq\\
  &=& <g \ot 1 \ot 1_\clq, (1 \ot \alpha(f))(k \ot 1_\clq)>_\clq\\
  &=& <g \ot 1_\clq, ({\rm id} \ot \tau_\mu \ot {\rm id}_\clq)((1 \ot \alpha(f))(k \ot 1_\clq))>_\clq\\
  &=& <g \ot 1_\clq, (T \ot {\rm id})(\alpha(f))>_\clq\\
  &=& <g \ot 1_\clq, (T \ot {\rm id})(Uf)>_\clq.\\
  \eean
  This proves the claim. 
  
  The rest of the proof of the theorem will be very similar to that of Theorem 2.14 of \cite{goswami}, replacing the Laplacian by the operator $T$. 
Indeed, it is easy to see from the arguments of \cite{goswami} that there is a universal CQG, say $\clq_1$, which acts on $C(X)$ faithfully in a way that the action preserves $\tau_\mu$, and the induced unitary representation on $L^2(\mu)$, say $V$, commutes with $T$ in the sense that $VT=(T \ot {\rm Id})  \circ V $. Clearly, both $\clq$ and $C(G)$ are quantum subgroups of  $\clq_1$. Moreover, reversing the steps in the proof of (\ref{***}) above, we get (\ref{1972}) with $U$ replaced by $V$, and since $d=\phi^{-1}(k)$, we conclude that $Q_1$ is an object of the category ${\bf Q}_{X,d}$, so it is indeed the desired universal object of the category ${\bf Q}_{X,d}$.
  \qed
\brmrk
The hypotheses of Theorem \ref{exist} are valid for any finite group and typical compact Riemannian Lie groups  (with the unique invariant metric)   for which the heat kernel $k$ is a bijective function of the metric $d$ ( see \cite{ter} for examples of such manifolds which include compact simply connected semisimple Lie groups). In such situations, one can choose $\mu$ to be the Haar measure and 
 $T$ to be the heat semigroup $T_t$ for some $t>0$.
\ermrk  
  
 \bcrlre
 
 For $X=\IT^n$, $n \geq 1$,  the universal object $QISO(X,d)$ coincides with $C(ISO(X))$.

 \ecrlre 
{\it Proof:}\\ It suffices to observe that  one can choose $T$ to be the associated heat semigroup $T_t$, for any $t>0$. Then it follows from the construction of   $QISO(X,d)$ that it must be a quantum subgroup of the quantum isometry group corresponding to the Riemannian Laplacian in the sense of \cite{goswami}, and this latter quantum group is known to be the same as the classical one from 
\cite{jyotish} and \cite{bhowmick}.
\qed

{\bf Acknowledgment:}\\
The author would like to thank Adam Skalski, Shuzhou Wang and Jyotishman Bhowmick for many useful conversations and discussion. He is  grateful to  Adam and Shuzhou for reading an earlier draft of the paper with great care and patience, and 
 pointing out several gaps in proofs and typos. He also thanks an anonymous referee who noticed two serious mistakes in an earlier version.


\begin{thebibliography} {ter}

\bibitem{ter} Arede, T.: Manifolds for which the heat kernel is given in terms of geodesic lengths, \emph{Lett. Math. Phys.} {\bf 9} (1985), 121--131.\\

\bibitem{ban1} Banica, T.: Quantum automorphism groups of small metric spaces, Pacific J.\ Math.\ {\bf 219}(2005), no. 1, 27--51.\\

\bibitem{ban2} Banica, T.: Quantum automorphism groups of homogeneous graphs, J.\ Funct.\ Anal.\ {\bf 224}(2005), no. 2, 243--280.\\

\bibitem{bichon} Bichon, J.: Quantum automorphism groups of finite graphs,  Proc.\ Amer.\ Math.\ Soc.\ {\bf 131}(2003), no. 3, 665--673.\\ 

\bibitem{bichon_alg} Bichon, J.: Algebraic quantum permutation groups,  Asian-Eur. J.\  Math.\ {\bf 1}(2008),  no. 1, 1--13.\\

\bibitem{vaes} Bichon, J., De Rijdt, A.,  Vaes, S.:  Ergodic coactions with large multiplicity and monoidal equivalence of quantum groups,  Comm. Math. Phys.  {\bf 262 } (2006),  no. 3, 703--728.\\

\bibitem{bhowmick} Bhowmick, J.: Quantum isometry group of the $n$-tori, Proc. Amer. Math. Soc. {\bf 137 }(2009), 3155-3161.\\

\bibitem{erg} Boca, F.: Ergodic actions of compact matrix pseudogroups on $C^*$-algebras, Recent advances in operator algebras 
(Orléans, 1992);
Astérisque No. 232 (1995), 93--109.\\ 

\bibitem{goswami} Goswami, D.: Quantum Group of isometries in Classical and Non Commutative Geometry. Comm. Math.Phys.{\bf 285}(2009), no.1, 141-160.\\

\bibitem{jyotish} Bhowmick, J. and Goswami, D.: Quantum isometry groups : examples and computations,  Comm. Math. Phys. {\bf 285}(2009), no. 2, 421-444.\\



\bibitem{qdisc} Bhowmick, J., Goswami, D. and Skalski, A.: Quantum isometry groups of $0$-dimensional manifolds, to appear in Trans. A. M. S.
(arXiv 0807.4288).\\

\bibitem{qorient} Bhowmick, J. and Goswami, D.: Quantum group of orientation preserving Riemannian isometries, \emph{ J. Funct.  Anal. }{\bf 257} (2009), 2530--2572.\\

\bibitem {connes} Connes, A.: ``Noncommutative Geometry", Academic Press, London-New York (1994).\\




\bibitem{vnqgp} Kustermans, J.; Vaes, S.: Locally compact quantum groups in the von Neumann algebraic setting,   Math.\  Scand.\ {\bf 92}  (2003),  no. 1, 68--92.\\


\bibitem{podles} Podles, Piotr: Symmetries of quantum spaces. Subgroups and quotient spaces of quantum ${\rm SU}(2)$ and ${\rm SO}(3)$ groups,  Comm.\  Math.\  Phys.\ {\bf  170} (1995),  no. 1, 1--20.\\

\bibitem{sabbe} Sabbe, M. and Quaegebeur, J.:  Isometric coactions of compact quantum groups on compact quantum metric spaces, arXiv 1007.0363.\\


\bibitem{soltan} Soltan, P. M.: Quantum families of maps and quantum semigroups on finite quantum spaces, preprint, arXiv:math/0610922.\\



\bibitem{vandaelenotes} Maes, A. and Van Daele, A.: Notes on compact quantum groups, Nieuw Arch.\ Wisk.\ (4){\bf 16}(1998), no. 1-2, 73--112.\\
\bibitem{free} Wang, S.: Free products of compact
quantum groups,  Comm.\ Math.\ Phys.\ {\bf 167} (1995), no. 3, 671--692.\\

\bibitem{wang1} Wang, S.: Quantum symmetry groups of finite spaces, Comm.\ Math.\ Phys.\ {\bf 195}(1998), 195--211.\\

\bibitem{wang2} Wang, S.: Structure and isomorphism classification of compact quantum groups $A_u(Q)$ and $B_u(Q)$, J.\ Operator Theory {\bf 48} (2002), 573--583.\\

\bibitem{woro1} Woronowicz, S. L.: Compact matrix pseudogroups, Comm.\ Math.\ Phys.\ {\bf 111}(1987), no. 4, 613--665.\\

\bibitem{woro2} Woronowicz, S. L.:  ``Compact quantum groups", pp. 845--884 in  Sym\'etries quantiques (Quantum symmetries)
 (Les Houches, 1995), edited by A. Connes et al., Elsevier, Amsterdam,
 1998.\\
  
 
 \bibitem{woro_pseudo}Woronowicz, S. L.: Pseudogroups, pseudospaces and Pontryagin duality, Proceedings of the International Conference on Mathematical Physics, Lausane (1979), Lecture Notes in Physics {\bf 116}, pp. 407-412.\\
 
\end{thebibliography}
\end{document}